\documentclass[11pt]{amsart}
\usepackage{etex}
\usepackage{diagbox}
\usepackage{array}
\usepackage{graphicx, pictex}
\usepackage{stmaryrd}
\usepackage{float}
\restylefloat{table}
\usepackage{multirow}
\usepackage{amsmath,amssymb,amsthm,mathrsfs}
\allowdisplaybreaks
\usepackage[dvipsnames]{xcolor}
\def\d{\Omega}

\def\du#1#2#3{\overset{#3}{\underset{#2}{#1}}}
\def\Forall{\quad \hbox{ for all }}

\def\M{{\mathcal{M}}}
\def\T{\mathcal{T}}

\newcommand{\tn}[1]{\lVert\kern-1pt\lvert{#1}\rvert\kern-1pt\rVert}
\def\<{{\langle}}
\def\>{{\rangle}}
\def\Forall{\quad \hbox{ for all }}

\def\bn{{\mathbf n}}

\def\C{{\mathcal{C}}}

\def\d{\Omega}

\def\Forall{\quad \hbox{ for all }}

\def\d{\Omega}

\def\Forall{\quad \hbox{ for all }}

\def\tb#1{{\|\kern-1pt| #1 \|\kern-1pt|}}
\def\nm2#1#2{\|#1\|_{2,\d_{#2}}}

\def\R{\mathbb{R}}

\def\T{\mathcal{T}}

 \theoremstyle{plain}
 \newtheorem{thm}{Theorem}[section]
 \numberwithin{equation}{section} 
 \numberwithin{figure}{section} 
 \theoremstyle{plain}
 \newtheorem{prop}[thm]{Proposition}
 \theoremstyle{plain}
 \newtheorem{algorithm}[thm]{Algorithm} 
 \theoremstyle{plain}
 \newtheorem{theorem}[thm]{Theorem}
 \theoremstyle{plain}
 
\theoremstyle{plain}
 \newtheorem{remark}[thm]{Remark}
 \theoremstyle{plain}
 \newtheorem{lemma}[thm]{Lemma}

\def\C{{\mathcal{C}}}
\def\T{{\mathcal{T}}}
\def\M{{\mathcal{M}}}

\def\d{{\Omega}}

\def\Forall{\quad \hbox{ for all }}

\def\<{{\langle}}
\def\>{{\rangle}}
\def\R{\mathbb{R}}

\def\bn{{\mathbf n}}

\def\du#1#2#3{\overset{#3}{\underset{#2}{#1}}}

\usepackage{amssymb}

\begin{document}

\title[Saddle Point Least Squares for second order elliptic interface]
{A nonconforming saddle point least squares approach for elliptic interface problems}

\author{Constantin Bacuta}
\address{University of Delaware,
Department of Mathematics,
501 Ewing Hall 19716}
\email{bacuta@udel.edu}

\author{Jacob Jacavage}
\address{University of Delaware,
Department of Mathematics,
501 Ewing Hall 19716}
\email{jjacav@udel.edu}


\thanks{The work was supported  by NSF, DMS-1522454 \\ Note: The original manuscript was submitted to {\it De Gruyter, Computational Methods in Applied Mathematics} on December 19, 2017}

\begin{abstract}
We present a non-conforming least squares method for  approximating  solutions of second order elliptic problems with discontinuous coefficients. The method is based on a general Saddle Point Least Squares (SPLS) method introduced in previous work based on conforming discrete spaces. 
The SPLS method has the advantage that a discrete $\inf-\sup$ condition is automatically satisfied for standard choices of test and trial spaces. We explore the SPLS method for non-conforming finite element trial spaces which allow higher order approximation of the fluxes. For the proposed iterative solvers, inversion at each step requires bases only for the test spaces. We focus on using projection trial spaces  with local projections that are easy to compute. The choice of the  local projections for the trial space  can be combined with classical gradient recovery techniques to lead to quasi-optimal  approximations of the global flux. Numerical results for 2D and 3D domains are included to support the proposed method. 
\end{abstract}
\maketitle

\section{Introduction}
Elliptic interface problems have applications in a variety of different fields. In material science, they arise in the study and design of composite materials built from essentially different components, see \cite{babuska94, Gueribiz2013, widlund2003, bakhalov89}. In fluid dynamics, they model several layers of fluids with different viscosities or diffusion through heterogeneous porous media \cite{Bernardi-Verfurth2000, Efendiev2007}. In addition, the elliptic interface problem is used to model stationary heat conduction problems with a conduction coefficient which is discontinuous across a smooth internal interface \cite{hansbo2002}, as well as in biological systems \cite{interfacebio}.

Given $f\in L^2(\Omega)$, we consider the problem of finding $u \in H_0^1(\Omega)$ such that 
\begin{equation}\label{interface_prob}
-\text{div}(A\nabla u)=f \ \ \text{in} \ \Omega,
\end{equation}
where the matrix $A$ is uniformly coercive and the entries could be discontinuous across an interface contained in $\Omega$, with possibly large jumps, across the subdomain boundaries. We also assume the continuity of the co-normal derivative along the interface(s), see Section \ref{sec:EllipticInterface}.  

The primal mixed variational formulation of \eqref{interface_prob} we consider is: Find $p =A\nabla u$, with $u\in H_0^1(\Omega)$, such that
\begin{equation}\label{int_varform}
(p,\nabla v)=(A\nabla u,\nabla v)=(f,v) \Forall v\in H_0^1(\Omega).
\end{equation}
Introducing the auxiliary variable $w \in H_0^1(\Omega)$, a saddle point reformulation of \eqref{int_varform} is to find $(w=0, p) \in H_0^1(\Omega) \times A \nabla H_0^1(\Omega)$ such that 
\begin{equation}\label{eq:PD4interface}
\begin{array}{lclll}
(\nabla w,\nabla v) & + & (p,\nabla v) &= (f,v) &\ \Forall  v \in  H_0^1(\Omega),\\
A\nabla w & & & =0.
\end{array}
\end{equation}
By defining the spaces $V:=H_0^1(\Omega), Q:=A\nabla V$, and $\tilde{Q}:=L^2(\Omega)^d$, as well as defining the bilinear form $b:V\times \tilde{Q}\to \mathbb{R}$ by
\[
b(v,q):=(q,\nabla v) \Forall v\in V,q\in \tilde{Q}, 
\]
problem \eqref{int_varform} can be rewritten as: Find $p\in Q$ such that 
\begin{equation}\label{int_varform2}
b(v, p)=(f,v) \Forall v\in V.
\end{equation}
Furthermore, by denoting $a_0(u,v):=(\nabla u, \nabla v)$ as the standard inner product on $V$, the saddle point reformulation \eqref{eq:PD4interface} can be rewritten as: Find $(w=0, p) \in V\times Q$ such that 
\begin{equation}\label{eq:PD4interface2}
\begin{array}{lclll}
a_0(w, v) & + & b(v,p) &= (f,v) &\ \Forall  v \in  V,\\
b(w,q) & & & =0 &\ \Forall  q \in  Q.
\end{array}
\end{equation}

The advantage of reformulating  \eqref{int_varform} (or  \eqref{int_varform2})  into \eqref{eq:PD4interface} (or  \eqref{eq:PD4interface2}) resides in the fact that we can use non-conformoing discrete finite element spaces to approximate $p=A\nabla u$, which lead to a better approximation for $p$ if compared with a direct approximation for $u$ from a variational formulation of \eqref{interface_prob} followed by the application of the linear operator $A\nabla\cdot$. In addition, we can apply the classical approximation theory for saddle point problems.

This idea can be extended to a more general class of mixed variational problems, and in \cite{BQ15} it was called the {\it Saddle Point Least Squares} (SPLS) method. The version we propose in this paper can be applied to the interface problem \eqref{interface_prob}, as well as more general first or second order elliptic PDEs. The SPLS method {\it bridges} between the field of {\it least squares methods} and the field of {\it symmetric saddle point problems}. The discretization approach in this paper can be viewed as a new discontinuous Petrov-Galerkin method. From the point of view of choosing the discrete spaces, it can be characterized as a dual  of Demkowicz-Gopalakrishnan's Discontinuous Petrov-Galerkin (DPG) method \cite{demkowicz-gopalakrishnanDPG10, demkowicz-gopalakrishnanDPG13}, which is currently undergoing an intensive study. 

While both methods have strong connections with least squares and minimum residual techniques, our proposed discretization process stands apart from the DPG approach due to the different ways in which the trial and test spaces are chosen. In our approach, we choose a discrete test space first and the trial space is then built in order to satisfy a discrete $\inf-\sup$  condition. For the SPLS method, the trial space is built from the action of the continuous differential operator associated with the problem on the test space. Due to the iterative process we choose to solve the discrete SPLS formulation, assembly of the stiffness matrices for the trial spaces is avoided. The SPLS method can be also be combined with multilevel preconditioning techniques in order to address particular challenges of the PDE to be solved due to discontinuous coefficients or multidimensional domains \cite{BJprec}. In contrast with the SPLS work presented in \cite{BQ15,BQ17}, where both the test and trial spaces were chosen to be conforming finite element spaces, this paper considers trial spaces which are non-conforming finite element spaces. This allows efficient treatment of PDEs with discontinuous coefficients.  

The paper is organized as follows. In section 2, we introduce notation for the general non-conforming (n-c) SPLS method and present two types of trial spaces along with stability and approximability properties. In section 3, the general theory will be applied to approximating the solution of second order elliptic problems with discontinuous coefficients. In section 4, numerical results for the SPLS dicretization are presented.  

\section{The general non-conforming SPLS approach}\label{section:ReviewLSPP}
We first introduce some notation for the spaces and operators for the general abstract setting. Let $V$ and $\tilde{Q}$ be infinite dimensional Hilbert spaces and assume the inner products $a_0(\cdot, \cdot)$ and $(\cdot, \cdot)_{\tilde{Q}}$ induce the  norms $|\cdot|_V =|\cdot| =a_0(\cdot, \cdot)^{1/2}$ and $\|\cdot\|_{\tilde{Q}}=\|\cdot\|=(\cdot, \cdot)_{\tilde{Q}}^{1/2}$. We denote the duals of $V$ and $\tilde{Q}$ by $V^*$ and $\tilde{Q}^*$, respectively. The dual pairings on $V^* \times V$ and $\tilde{Q}^* \times \tilde{Q}$ will both be denoted by $\langle \cdot, \cdot \rangle$. With the inner product 
$(\cdot, \cdot)_{\tilde{Q}}$, we associate the operator $\C: \tilde{Q}\to\tilde{Q}^*$ defined by 
\[
\langle \C p,q\rangle=(p, q)_{\tilde{Q}} \Forall  \ p, q \in \tilde{Q}.
\]
The operator $\C^{-1}:\tilde{Q}^*\to \tilde{Q}$ is the Riesz-canonical isometry. In addition, we let $Q$ be a closed subspace of $\tilde{Q}$ equipped with the induced inner product (from $\tilde{Q}$).

We assume that $b(\cdot, \cdot)$ is a continuous bilinear form on $V\times \tilde{Q}$ satisfying
\begin{equation}\label{sup-sup_a}
\du{\sup}{p \in \tilde{Q}}{} \ \du {\sup} {v \in V}{} \ \frac {b(v, p)}{|v|\,\|p\|} =M <\infty, 
\end{equation}
and the following $\inf-\sup$ condition on $V\times Q$,
 \begin{equation} \label{inf-sup_a}
\du{\inf}{p \in Q}{} \ \du {\sup} {v \in V}{} \ \frac {b(v, p)}{|v|\,\|p\|} =m>0.
\end{equation} 
With the form $b$, we associate the linear operators $B:V\to \tilde{Q}^*$ and $B^*:\tilde{Q}\to V^*$ defined by
\[
\<B v,q\>=b(v, q)= \langle B^*q, v \rangle  \Forall v \in V, \ q \in \tilde{Q}.
\]
Lastly, we define $V_0$ to be the kernel of $B$, i.e.,
\[ 
V_0 :=\mathrm{Ker}(B)= \{v \in V |\  Bv=0\}.
\]

We consider problems of the form: Given $F \in V^*$, find $p \in Q$ such that
\begin{equation}\label{eq:BBprobb}
b(v, p) = \<F,v\>  \Forall v \in V, \ \ \ \ \text{or}\ \ \  \ B^*p=F. 
\end{equation}
We note here that for the existence and uniqueness of the solution of the continuous problem \eqref{eq:BBprobb}, we use the trial space $Q$. However, for discretization purposes, we need to consider the form $b(\cdot, \cdot)$ on $V \times \tilde{Q}$.
The existence and uniqueness of \eqref{eq:BBprobb} was first studied by Aziz and Babu\v{s}ka in \cite{A-B}.  It is well known that if a bounded form $b:V \times \tilde{Q} \to \R$ satisfies \eqref{inf-sup_a} and the data $F \in V^*$ satisfies the {\it compatibility condition} 
\begin{equation}\label{eq:BBsuf}
\<F,v\> =0 \Forall v \in V_0, 
\end{equation}
then the mixed problem \eqref{eq:BBprobb} has a unique solution, see e.g. \cite{A-B,B09}. With the mixed problem  \eqref{eq:BBprobb}, we associate the SPLS formulation: Find $(w, p) \in (V, Q)$ such that  
\begin{equation}\label{abstract:variational2}
\begin{array}{lclll}
a_0(w,v) & + & b( v, p) &= \langle F,v \rangle &\ \Forall  v \in V,\\
b(w,q) & & & =0   &\  \Forall  q \in Q.  
\end{array}
\end{equation}
The following statement summarizes the connection between the two variational formulations. The remark was pointed out  in \cite{BM12, Dahmen-Welper-Cohen12} and  is essential in our approach and (some versions of) the DPG method. It is worth noting that the $p$ component of the solution of \eqref{abstract:variational2} is in fact the solution of the normal equation that corresponds to our main problem  \eqref{eq:BBprobb}, see \cite{BQ15}. 

\begin{prop}\label{prop:SPLS}
 In the presence of the continuous $\inf-\sup$ condition \eqref{inf-sup_a} and the {compatibility condition} \eqref{eq:BBsuf}, we have that $p$ is the unique solution of \eqref{eq:BBprobb}  {if and only if} $(w=0 , p) $ is the unique solution of \eqref{abstract:variational2}. 
\end{prop} 

\subsection{Non-Conforming SPLS discretization} \label{sec:SPLSd}

The non-conforming (trial space) {\it SPLS  discretization} of \eqref{eq:BBprobb} is defined as a (trial) non-conforming saddle point discretization of \eqref{abstract:variational2}.  We consider finite dimensional approximation spaces $V_h \subset V$ and $\M_h \subset \tilde{Q}$ (larger  than $Q$ in general) and  restrict  the forms $a_0(\cdot, \cdot)$ and $b(\cdot, \cdot)$ to the discrete spaces $V_h$ and $\M_h$. Assume that the following discrete $\inf-\sup$ condition holds for the pair $(V_h, \M_h)$: 
 \begin{equation}
\label{inf-sup_h}
\du{\inf}{p_h \in \M_h}{} \ \du {\sup} {v_h \in V_h}{} \ \frac {b(v_h, p_h)}{|v_h|\,\|p_h\|} =m_h>0.
\end{equation} 
 We define $V_{h,0}$ to be the kernel of the discrete operator $B_h$, i.e., 
 \[
 V_{h,0}:= \{v_h \in V_h | \ b(v_h, q_h) =0 \Forall q_h  \in\M_h\},
 \]
and let $V_{h,0}^\perp$ denote the orthogonal complement of $V_{h,0}$ with respect to the inner product $a_0(\cdot, \cdot)$ on $V_h$.  If $V_{h,0} \subset V_0$, then the compatibility condition \eqref{eq:BBsuf} implies a discrete compatibility condition. Consequently, under the discrete stability assumption \eqref{inf-sup_h}, the problem of finding $p_h \in \M_h$ such that  
\begin{equation}\label{eq:BBprobb-hh}
b(v_h, p_h) = \<F, v_h\> \Forall v_h \in V_h,
\end{equation}
has a unique solution.

In general, the compatibility condition \eqref{eq:BBsuf} might not hold on $V_{h,0}$. Hence, the discrete problem \eqref{eq:BBprobb-hh} may not be well-posed. In any case, under the assumption \eqref{inf-sup_h}, the standard discrete saddle point problem of finding $(w_h, p_h) \in V_h \times \M_h$ such that  
\begin{equation}\label{abstract:variational2_h}
\begin{array}{lclll}
a_0(w_h,v_h) & + & b( v_h, p_h) &= \langle F,v_h \rangle &\ \Forall  v_h \in V_h,\\
b(w_h,q_h) & & & =0   &\  \Forall  q_h \in\M_h, 
\end{array}
\end{equation}
does have a unique solution. We call the variational formulation \eqref{abstract:variational2_h} the {\it non-conforming saddle point least squares discretrization} of \eqref{eq:BBprobb}. As in the continuous case, it is easy to prove that the $p_h$ part of the solution of \eqref{abstract:variational2_h} is the solution of the {\it normal equation} associated with \eqref{eq:BBprobb-hh}. 

\subsection{The discrete spaces}\label{sec:Mh}
Let $V_h$ be a {\it finite element subspace} of $V$ and assume that the action of $\C^{-1}$ at the continuous level is easy to obtain. 

\subsubsection{No projection trial space}\label{sec:noproj}
We first consider the case when $\M_h$ is given by
\[
\M_h:=\C^{-1}BV_h \subset  \tilde{Q}.
\]
In this case, we have $V_{h,0}\subset V_0$ and a discrete $\inf-\sup$ condition holds. Indeed, for a generic $p_h =\C^{-1}Bw_h\in \M_h$ where $w_h\in V_{h,0}^{\perp}$, we have 
\begin{align*}
\inf_{p_h \in \M_h}\sup_{v_h\in V_h} \frac{b(v_h,p_h)}{|v_h|\  \|p_h\|}&=\inf_{w_h \in V^{\perp}_{h,0}}\sup_{v_h\in V_h} \frac{(\C^{-1}Bv_h,\C^{-1}Bw_h)_{\tilde{Q}}}{|v_h|\,\|\C^{-1}Bw_h\|}\\
&\geq \inf_{w_h \in V^{\perp}_{h,0}}\frac{\|\C^{-1}Bw_h\|^2}{|w_h|\,\|\C^{-1}Bw_h\|} 
=\inf_{w_h \in V^{\perp}_{h,0}}\frac{\|\C^{-1}Bw_h\|}{|w_h|}
:=m_{h,0}. \stepcounter{equation}\tag{\theequation}\label{inf-supNoProj}
\end{align*}
Thus, we have that both  variational formulations \eqref{eq:BBprobb-hh} and  \eqref{abstract:variational2_h} have a unique solution $p_h \in \M_h$. Furthermore, using Proposition \ref{prop:SPLS} for the discrete pair $(V_h, \M_h)$, we have that $(w_h =0, p_h) $ is the solution of \eqref{abstract:variational2_h}.

\subsubsection{Approximability of no projection trial space}\label{sec:noprojApprox}
Note that if $p$ is the solution of \eqref{eq:BBprobb} and $p_h$ is the solution of \eqref{eq:BBprobb-hh}, or $(0, p_h) $ is  the solution of \eqref{abstract:variational2_h}, then from \eqref{eq:BBprobb} and \eqref{eq:BBprobb-hh} we obtain
\[
0=b(v_h,p-p_h)=(\C^{-1}Bv_h,p-p_h)_{\tilde{Q}} \Forall v_h \in V_h.
\]
Thus, $p_h$ is the orthogonal projection of $p$ onto $\M_h$ which gives us
\[
\|p-p_h\|=\inf_{q_h\in \M_h}\|p-q_h\|.
\]
This result is optimal, and in contrast with the standard approximation estimates for saddle point problems, it does not depend on $m_{h,0}$.

\subsubsection{Projection type trial space}\label{sec:proj}
Let $\tilde{\M}_h \subset \tilde{Q}$ be a finite dimensional subspace equipped with the inner product $(\cdot,\cdot)_h$. Define the representation operator $R_h:\tilde{Q}\to \tilde{\M}_h$ by 
\begin{equation}\label{rep_operator}
(R_hp,q_h)_h:=(p,q_h)_{\tilde{Q}} \Forall q_h \in \tilde{\M}_h.
\end{equation}
Here, $R_hp$ is the Riesz representation of $p\to(p,q_h)_{\tilde{Q}}$ as a functional on $(\tilde{\M}_h,(\cdot,\cdot)_h)$. 
\begin{remark}
In the case when $(\cdot,\cdot)_h$ coincides with the inner product on $\tilde{Q}$, we have that $R_h$ is the orthogonal projection onto $\tilde{\M}_h$.
\end{remark}

Since the space $\tilde{\M}_h$ is finite dimensional, there exist constants $k_1,k_2$ such that 
\begin{equation}\label{norm-equiv}
k_1\|q_h\|\leq \|q_h\|_h\leq k_2\|q_h\| \Forall q_h \in \tilde{\M}_h.
\end{equation}
We further assume that the equivalence is uniform with respect to $h$, i.e., the constants $k_1, k_2$ are independent of $h$. 
Using the operator $R_h$, we define $\M_h$ as
\[
\M_h:=R_h\C^{-1}BV_h\subset \tilde{\M}_h \subset \tilde{Q}.
\]
The following proposition gives a sufficient condition on $R_h$ to ensure the discrete $\inf-\sup$ condition is satisfied and relates the stability of the families of spaces $\{(V_h, \C^{-1} B V_h)\}$ and $\{(V_h, R_h \C^{-1} B V_h)\}$.
\begin{prop}\label{discrete-Stability}
Assume that 
\begin{equation}\label{eq:ProjCoerces}
\|R_h  q_h \|_h \geq \tilde{c}\,  \| q_h\| \Forall q_h \in \C^{-1} B V_h,
\end{equation}
with a constant $\tilde{c}$ independent of $h$. Then $V_{h,0} \subset V_0$. Furthermore, the stability of the family $\{(V_h, \C^{-1} B V_h) \}$, meaning $m_{h,0}$ defined in \eqref{inf-supNoProj} satisfies $m_{h,0} >c_0>0$ for some constant $c_0$ independent of $h$, implies the stability of the family $ \{(V_h, R_h \C^{-1} B V_h)\}$.
\begin{proof}
Let $v_h \in V_{h,0}$. Then, for any $p_h \in \M_h$,  
\[
0=b(v_h,p_h)=(\C^{-1}Bv_h,p_h)_{\tilde{Q}}=(R_h\C^{-1}Bv_h,p_h)_h.
\]
Taking $p_h=R_h\C^{-1}Bv_h$ gives us $\|R_h\C^{-1}Bv_h\|_h=0$ and the inclusion $V_{h,0} \subset V_0$ follows from \eqref{eq:ProjCoerces}. To show the stability, we take a generic function $p_h =R_h\C^{-1}Bw_h\in \M_h$ where $w_h \in V^{\perp}_{h,0}$. We have
\begin{align*}
m_h=\inf_{p_h \in \M_h}\sup_{v_h\in V_h} \frac{b(v_h,p_h)}{|v_h|\  \|p_h\|_h}&= \inf_{w_h \in V^{\perp}_{h,0}}\sup_{v_h\in V_h} \frac{(\C^{-1}Bv_h,R_h\C^{-1}Bw_h)_{\tilde{Q}}}{|v_h|\  \|R_h\C^{-1}Bw_h\|_h}\\
&=\inf_{w_h \in V^{\perp}_{h,0}}\sup_{v_h\in V_h} \frac{(R_h\C^{-1}Bv_h,R_h\C^{-1}Bw_h)_h}{|v_h|\  \|R_h\C^{-1}Bw_h\|_h}\\
&\geq \inf_{w_h \in V^{\perp}_{h,0}} \frac{\|R_h\C^{-1}Bw_h\|_h^2}{|w_h|\  \|R_h\C^{-1}Bw_h\|_h}\\
&\geq \tilde{c} \inf_{w_h \in V^{\perp}_{h,0}} \frac{\|\C^{-1}Bw_h\|}{|w_h|}
=\tilde{c}\, m_{h,0},
\end{align*}
where $m_{h,0}$ is defined in \eqref{inf-supNoProj}.
\end{proof}
\end{prop}
As a consequence of Proposition \ref{discrete-Stability}, we have that under the assumption \eqref{eq:ProjCoerces} both variational formulations \eqref{eq:BBprobb-hh} and  \eqref{abstract:variational2_h} have  unique solution $p_h \in \M_h$. Furthermore, using Proposition \ref{prop:SPLS} for the discrete pair $(V_h, \M_h)$, we have that $(w_h =0, p_h) $ is the solution of \eqref{abstract:variational2_h}.

\subsubsection{Approximability of projection type trial space} The following proposition shows that under condition \eqref{eq:ProjCoerces} we have a quasi-optimal approximability property for the projection type trial space. 
\begin{prop}\label{approxProj}
If $p$ is the solution of \eqref{eq:BBprobb}, $p_h$ is the solution of \eqref{eq:BBprobb-hh} (or the n-c SPLS solution of \eqref{abstract:variational2_h}), and $R_h$ satisfies \eqref{eq:ProjCoerces}, then
\[
\|p-p_h\| \leq C\inf_{q_h\in \M_h}\|p-q_h\|,
\]
where $C$ depends only on $\tilde{c}$ of \eqref{eq:ProjCoerces} and the equivalence of norms constants of \eqref{norm-equiv}.
\begin{proof}
From the assumptions on $p$ and $p_h$, using \eqref{eq:BBprobb} and \eqref{eq:BBprobb-hh} we obtain
\[
0=b(v_h, p-p_h) = (\C^{-1} Bv_h, p-p_h)_{\tilde{Q}} \Forall \ v_h \in V_h. 
\]
In turn, this implies
\begin{equation}\label{Q_innerprod_equality}
(\C^{-1}Bv_h,p-Q_hp)_{\tilde{Q}}=(\C^{-1}Bv_h,p_h-Q_hp)_{\tilde{Q}} \Forall v_h\in V_h,
\end{equation}
where $Q_h$ is the orthogonal projection onto $\M_h$. Note that
\begin{equation}\label{HilbertSpaceNorm}
\|p_h-Q_hp\|_h=\sup_{q_h\in \M_h}\frac {|(p_h-Q_hp,q_h)_h|}{\|q_h\|_h}.
\end{equation}
Using \eqref{eq:ProjCoerces} and \eqref{Q_innerprod_equality}, we obtain
\begin{align*}
\sup_{q_h\in \M_h}\frac {|(p_h-Q_hp,q_h)_h|}{\|q_h\|_h}&=\sup_{w_h\in V^{\perp}_{h,0}}\frac {|(p_h-Q_hp,R_h\C^{-1}Bw_h)_h|}{\|R_h\C^{-1}Bw_h\|_h}\\
&=\sup_{w_h\in V^{\perp}_{h,0}}\frac {|(p_h-Q_hp,\C^{-1}Bw_h)_{\tilde{Q}}|}{\|R_h\C^{-1}Bw_h\|_h}\\
&=\sup_{w_h\in V^{\perp}_{h,0}}\frac {|(p-Q_hp,\C^{-1}Bw_h)_{\tilde{Q}}|}{\|R_h\C^{-1}Bw_h\|_h}\\
&\leq\sup_{w_h\in V^{\perp}_{h,0}}\frac {\|p-Q_hp\| \,\|\C^{-1}Bw_h\|}{\|R_h\C^{-1}Bw_h\|_h} 
\leq \frac{1}{\tilde{c}}\|p-Q_hp\|.
\end{align*}
Hence, from \eqref{norm-equiv}, \eqref{HilbertSpaceNorm}, and the above estimate we have
\begin{equation}\label{norm_equiv_estimate}
\|Q_hp-p_h\|\leq \frac{1}{k_1}\|Q_hp-p_h\|_h\leq \frac{1}{\tilde{c}k_1}\|p-Q_hp\|.
\end{equation}
Thus, 
\begin{align*}
\|p-p_h\| &\leq \|p-Q_hp\|+\|Q_hp-p_h\|\\
&\leq \left(1+\frac{1}{\tilde{c}k_1}\right) \|p-Q_hp\|
=C\, \inf_{q_h\in \M_h}\|p-q_h\|.
\end{align*}
\end{proof}
\end{prop}

\begin{remark}
The no projection trial space described in Section \ref{sec:noproj} can be viewed as the special case of the projection type trial space when $R_h=I$.
\end{remark}

\subsection{Iterative solvers}\label{sec:solvers}
When solving \eqref{abstract:variational2_h} on $(V_h,\M_h=\C^{-1}BV_h)$ or $(V_h,\M_h=R_h\C^{-1}BV_h)$, a global linear system might be difficult to assemble as one may not be able to find simple local bases for the space $\M_h$, especially for the projection type of trial space. Nevertheless, it is possible to solve \eqref{abstract:variational2_h} without an explicit basis for $\M_h$  choice of by using the Uzawa (U), Uzawa Gradient (UG), or Uzawa Conjugate Gradient (UCG) algorithm. We will describe each algorithm below. For implementation and  convergence analysis for such algorithms, it is essential to use the $(\cdot,\cdot)_h$ inner product on $\M_h$.

\begin{algorithm} (U-UG) Algorithms \label{alg:U-UG}
\vspace{0.1in}

{\bf Step 1:} {\bf Set} $p_0=0 \in \M_h$, 
{\bf compute} $w_{1} \in V_h $, $q_1 \in \M_h$ by  
\[
\begin{aligned}
& a_0( w_{1}, v) &= &\ \langle F,v\rangle - b(v, p_{0}) &\ &\Forall v \in V_h,\\
&  (q_1, q)_{h} & =  &\ b(w_1 ,q)  &\ &\Forall  q \in \M_h.
\end{aligned}
\]
\vspace{0.1in}

{\bf Step 2:} {\bf For} $j=1,2,\ldots, $ {\bf compute} $h_j, \alpha_j, p_j, w_{j+1}, q_{j+1}$ by 
\[
\begin{aligned}
& {\bf (U-UG1)} \ \ \ \  & a_0( h_{j}, v) &= - b(v, q_j)  && \text{for all} \ v \in V_h\\
& {\bf (U\alpha)} \ \ \ \  & \alpha_j &=\ \alpha_0  && \text{for the Uzawa algorithm or}\\
& {\bf (UG\alpha)} \ \ \ \  & \alpha_j &= - \frac{(q_j, q_j)_{h}}{b(h_j,q_j)}  && \text{for the UG algorithm} \\
& {\bf (U-UG2)} \ \ \ \  & p_{j} &= \ p_{j-1} + \alpha_j \  q_j \\
& {\bf (U-UG3)} \ \ \ \  & w_{j+1} &= \ w_j + \alpha_j\ h_j \\
& {\bf (U-UG4)} \ \ \ \  & (q_{j+1}, q)_{h} &=  \ b(w_{j+1} ,q)   && \text{for all}  \ q \in \M_h.
\end{aligned}
\]
\end{algorithm}
To obtain the UCG algorithm, the UG algorithm is modified as in \cite{braess, verfurtth84} by the following steps.
First, we define $d_1:=q_1$ in {\bf Step 1}. Then, we modify  {\bf Step 2} by replacing $ b(\cdot , q_j)$ with  $b(\cdot, d_j)$, where $\{d_j\}$ is a sequence of conjugate directions. The resulting algorithm is as follows:

\begin{algorithm} (UCG) Algorithm \label{alg:UCG}
\vspace{0.1in}

{\bf Step 1:} {\bf Set} $p_0=0 \in \M_h$. {\bf Compute} $w_{1} \in V_h $, $q_1, d_1 \in \M_h$ by 
\[ 
 \begin{aligned}
& a_0( w_{1}, v)& = &\ \<F,v\> - b(v, p_{0}) &\ &\Forall v \in V_h,&\\
&   (q_1, q)_{h} & = & \ b(w_1 ,q)  &\ &\Forall  q \in \M_h,& \ \ d_1:=q_1.
\end{aligned}
\]
\vspace{0.1in}

{\bf Step 2:} {\bf For} $j=1,2,\ldots, $ {\bf compute} 
 $h_j, \alpha_j, p_j, w_{j+1}, q_{j+1}, \beta_j, d_{j+1}$ by 
\[ 
 \begin{aligned}
& {\bf (UCG1)} \ \ \ \  & a_0( h_{j}, v) = & - b(v, d_j)  && \text{for all} \ v \in V_h&\\
& {\bf (UCG\alpha)} \ \ \ \  & \alpha_j =& - \frac{(q_j, q_j)_{h}}{b(h_j,q_j)} \\
& {\bf (UCG2)} \ \ \ \  & p_{j} = & \ p_{j-1} + \alpha_j \  d_j \\
& {\bf (UCG3)} \ \ \ \  & w_{j+1} = & \ w_j + \alpha_j\ h_j \\
& {\bf (UCG4)} \ \ \ \  & (q_{j+1}, q)_{h} = & \ b(w_{j+1} ,q) && \text{for all} \ q \in \M_h& \\
& {\bf (UCG\beta)} \ \ \ \  & \beta_j=& \ \frac{(q_{j+1}, q_{j+1})_{h}}{(q_j,q_j)_{h}} \\
& {\bf (UCG6)} \ \ \ \  & d_{j+1}= & \ q_{j+1} +\beta_j d_j. \\
\end{aligned}
\]
\end{algorithm}
Note that at each iteration step, only one inversion involving the form $a_0(\cdot,\cdot)$ is required. In \cite{B14}, it was proven that if $(w_h,p_h)$ is the discrete solution of \eqref{abstract:variational2_h} and $(w_{j+1},p_j)$ is the $j^{th}$ iteration for the U, UG, or UCG algorithm, then $(w_{j+1}, p_j)\to (w_h, p_h)$. In addition, there are constants $c_1,c_2$, independent of $h$, such that  for all $j=1,2, \dots$, we have
\begin{equation}\label{alg_estimate}
\begin{aligned}
\frac{c_1}{M^2}\|q_{j+1}\|\leq \|p_j-p_h\|\leq \frac{c_2}{m_h^2}\|q_{j+1}\|, \\
c_1 \, \frac{m_h}{M^2}\|q_{j+1}\|\leq \|w_{j+1}-w_h\|\leq c_2\, \frac{M}{m_h^2}\|q_{j+1}\|.
\end{aligned}
\end{equation}
Hence, the first equation in \eqref{alg_estimate} entitles $\|q_{j+1}\|$ as a computable, efficient, and uniform iteration error estimator for all three algorithms.  

We note that for the no projection choice of trial space $\M_h$ outlined in Section \ref{sec:noproj}, the residual $q_{j+1}$ from {\bf Step 1}, {\bf (U-UG4)}, and {\bf UCG4} can be computed using the action of the operator $\C^{-1}B$, i.e.,
\[
q_{j+1}=\C^{-1}Bw_{j+1}.
\]
Also, for the choice of a projection type trial space for $\M_h$ outlined in Section \ref{sec:proj}, the residual $q_{j+1}$ can be computed by applying the operator $\C^{-1}B$ followed by the operator $R_h$, i.e., 
\[
q_{j+1}=R_h(\C^{-1}Bw_{j+1}).
\]

\begin{remark}\label{U-remark}
If we focus only on the updates for the $p_j$'s in the U, UG, or UCG algorithm, we can see that they correspond to applying the standard Uzawa, Uzawa Gradient, or Uzawa Conjugate Gradient Algorithms (respectively) for inverting the Schur complement $S_h$ corresponding to the discrete system \eqref{abstract:variational2_h}. Due to the assumption \eqref{inf-sup_h}, $S_h$ is a symmetric positive definite operator. Consequently, the iterations $p_{j}$ converge to the solution $p_h$  with a rate of convergence that depends on the condition number of $S_h$, which is $\kappa(S_h) \leq \frac{M^2}{m_h^2}$. 
\end{remark}

\section{n-c SPLS for second order elliptic interface problems} \label{sec:EllipticInterface}
Let $\Omega \subset \mathbb{R}^d$ be a bounded polygonal domain with $\{\Omega_j\}_{j=1}^N$ a partition of $\Omega$ and $\bn_j$ be the outward unit normal vector to $\partial \Omega_j$. Define $\Gamma_{km}:=\partial \Omega_k \cap \partial\Omega_m$ to be the interface between $\Omega_k$ and $\Omega_m$ for $1\leq k<m\leq N$. Given $f\in L^2(\Omega)$, we consider the problem of finding $u \in H_0^1(\Omega)$ such that \eqref{interface_prob} is satisfied together with 
the continuity of the co-normal derivative condition
\[
\llbracket A\nabla u \cdot \bn\rrbracket_{\Gamma_{km}}=(A_k\nabla u_k\cdot \bn_k+A_m\nabla u_m\cdot \bn_m)\big|_{\Gamma_{km}}=0 \Forall k<m.
\]
We assume the matrix $A$ is symmetric and satisfies
\begin{equation}\label{Acoeff}
a_{min}|\xi|_e^2\leq \langle A(x)\xi,\xi\rangle_e\leq a_{max}|\xi|_e^2 \Forall x\in \Omega,\xi \in \mathbb{R}^d,
\end{equation}
for positive constants $a_{min}\leq a_{max}$ and where $\langle \cdot,\cdot\rangle_e$ and $|\cdot|_e$ denote the standard Euclidean inner product and norm for vectors in $\mathbb{R}^d$. In addition, the entries could by discontinuous, with possibly large jumps, across the subdomain boundaries. Throughout this section, $(\cdot,\cdot)$ and $\|\cdot\|$ will denote the standard $L^2$ inner product and norm for both scalar and vector functions. The primal mixed variational formulation of \eqref{interface_prob} we consider is: Find $p =A\nabla u$, with $u\in H_0^1(\Omega)$, such that
\eqref{int_varform} holds. 
We note here that, from the ellipticity assumption for $A$ and the fact that $ \nabla \cdot$ is injective on $H_0^1(\Omega)$, the representation of $p=A \nabla u$ is unique.

To fit \eqref{int_varform} into the abstract formulation \eqref{eq:BBprobb}, we let $V:=H_0^1(\Omega), \tilde{Q}:=L^2(\Omega)^d, Q:=A\nabla V$, and define $b:V\times \tilde{Q}\to \mathbb{R}$ by 
\[
b(v,q):=(q,\nabla v) \Forall v\in V,q\in \tilde{Q}.
\]
Also, 
\[
\langle F,v\rangle:=(f,v) \Forall v\in V.
\]
On $V$, we consider the standard inner product
\[
a_0(u,v):=(\nabla u, \nabla v) \Forall u,v \in V,
\]
and on $\tilde{Q}$, we define the weighted inner product
\[
(p,q)_{\tilde{Q}}:=(p,A^{-1}q) \Forall p,q\in \tilde{Q}.
\]
Note that for $\tau_1,\tau_2\in Q$, we then have 
\[
(\tau_1,\tau_2)_Q=(\tau_1,\tau_2)_{\tilde{Q}}=(A\nabla u_1,A\nabla u_2)_{\tilde{Q}}=(A\nabla u_1,\nabla u_2).
\]
With these inner products on $V$ and $\tilde{Q}$, we have that the operators $B:V\to \tilde{Q}^*$ and $\C^{-1}B:V\to \tilde{Q}$ are given by
\[
Bv=\nabla v, \ \ \text{and} \ \ \C^{-1}Bv=A\nabla v \Forall v\in V.
\]
Hence,
\[
V_0=Ker(B)=\{v\in V|Bv=0\}=\{v\in H_0^1(\Omega)|\nabla v=0\}=\{0\},
\]
which implies \eqref{eq:BBsuf} is trivially satisfied. We note that, as presented in \cite{BQ17}, the continuity constant satisfies
\begin{align*}
M&=\sup_{q\in \tilde{Q}}\sup_{v\in V}\frac{b(v,q)}{|v|_V\,\|q\|_{\tilde{Q}}}=\sup_{q\in \tilde{Q}}\sup_{v\in V}\frac{(q,\nabla v)}{|v|_V\,\|q\|_{\tilde{Q}}}\\
&=\sup_{q\in \tilde{Q}}\sup_{v\in V}\frac{(q,A\nabla v)_{\tilde{Q}}}{|v|_V\,\|q\|_{\tilde{Q}}}\leq \sup_{v\in V}\frac{\|A\nabla v\|_{\tilde{Q}}}{\|\nabla v\|}\leq \sqrt{a_{max}}< \infty. \stepcounter{equation}\tag{\theequation}\label{supsup_int}
\end{align*}
and  the $\inf-\sup$ constant satisfies
\begin{equation}\label{infsup_int}
 \begin{aligned}
m &=\inf_{q=A\nabla u \in Q}\sup_{v\in V}\frac{b(v,q)}{|v|_V\,\|q\|_{\tilde{Q}}} 
= \du{\inf}{u \in V}{} \ \du {\sup} {v \in V}{} \ \frac {(A \nabla u, \nabla v)}{(A \nabla u, \nabla u)^{1/2} \ |v|_V} \\ 
& \geq  \du{\inf}{u \in V}{}  \frac {(A \nabla u, \nabla u)}{(A \nabla u, \nabla u)^{1/2} \ \|\nabla u\|} 
\geq \sqrt{a_{min}}>0.
 \end{aligned}
\end{equation}
Consequently, the variational problem \eqref{int_varform} is well-posed and suitable for n-c SPLS formulation and discretization.

\subsection{n-c SPLS discretization for second order elliptic interface problems} 
We take $V_h \subset V=H_0^1(\Omega)$ to be the space of continuous piecewise polynomials of degree $k$ with respect to the {\it interface-fitted} triangular mesh $\T_h$. We note that while the {\it no projection trial space} case is similar with the work presented in \cite{BQ17}, the {\it  projection trial space} is  analyzed using the non-conforming trial space setting and leads to new stability and approximability estimates for the discontinuous coefficients (or interface) case. 

\subsubsection{No projection trial space}
Following Section \ref{sec:noproj}, we define the trial space as
\[
\M_h:=\C^{-1}BV_h=A\nabla V_h.
\]
By similar arguments used to show \eqref{infsup_int}, we obtain
\begin{equation}\label{infsup_discrete}
m_h:=\inf_{q_h=A\nabla u_h\in \M_h}\sup_{v_h\in V_h}\frac{b(v_h,q_h)}{|v_h|_V\,\|q_h\|_{\tilde{Q}}}\geq \sqrt{a_{min}}>0.
\end{equation}
Thus, we do have stability in this case. The discrete mixed variational formulation is: Find $p_h=A\nabla u_h$, with $u_h\in V_h$, such that
\begin{equation}\label{discrete_problem_noproj}
(p_h,\nabla v_h)=(A\nabla u_h,\nabla v_h)=(f,v_h) \Forall v_h\in V_h.
\end{equation}
The SPLS discretization \eqref{abstract:variational2_h} to be solved is: Find $(w_h,p_h=A\nabla u_h)$ such that
\begin{equation}\label{discrete_problemSPLS_noproj}
\begin{array}{lclll}
(\nabla w_h,\nabla v_h) & + & (p_h,\nabla v_h) &= (f,v_h) &\ \Forall  v_h \in V_h,\\
A\nabla w_h & & & =0.
\end{array}
\end{equation}

\subsubsection{Projection type trial space}\label{projForInt}
We define $\tilde{\M}_h\subset  \tilde{Q} =L^2(\Omega)^d$ to be
\[
\tilde{\M}_h:=\bigoplus_{i=1}^N AM_{h,0}|_{\Omega_i},
\]
where $N$ is the number of subdomains and where each component of $M_{h,0}|_{\Omega_i}$ consists of continuous piecewise polynomials of degree $k$ with respect to the mesh $\T_{h,i}:=\T_h|_{\Omega_i}$ with no restrictions on the boundary. We equip $\tilde{\M}_h$ with the inner product
\[
(A\tilde{q}_h,A\tilde{p}_h)_h=\sum_{i=1}^N(A\tilde{q}_h,A\tilde{p}_h)_{\tilde{Q},\Omega_i} \Forall A\tilde{q}_h,A\tilde{p}_h \in \tilde{\M}_h.
\]
Here, $(\cdot,\cdot)_{\tilde{Q},\Omega_i}$ is the inner product on $\tilde{Q}$ restricted to the subdomain $\Omega_i$. Using the definition of $R_h$ given in \eqref{rep_operator}, we have that for $p\in \tilde{Q}$
\begin{align*}
(p,A\tilde{q}_h)_{\tilde{Q}}=(R_hp,A\tilde{q}_h)_h&=\sum_{i=1}^N(R_hp,A\tilde{q}_h)_{\tilde{Q},\Omega_i}\\
&=(R_hp,A\tilde{q}_h)_{\tilde{Q}} \Forall A\tilde{q}_h\in \tilde{\M_h}.
\end{align*}
Thus, $R_hp$ is the orthogonal projection of $p$ onto $\tilde{\M_h}$ in the $(\cdot,\cdot)_{\tilde{Q}}$ inner product. In turn, this implies $R_hp|_{\Omega_j}$ is the orthogonal projection onto $\tilde{\M}_h|_{\Omega_j}=AM_{h,0}|_{\Omega_j}$ in the $(\cdot,\cdot)_{\tilde{Q}}$ inner product. We then define
\[
\M_h:=R_hA\nabla V_h.
\]
The discrete mixed variational formulation in this case is: Find $p_h=R_hA\nabla u_h$, with $u_h\in V_h$, such that
\begin{equation}\label{discrete_problem}
(p_h,\nabla v_h)=(R_hA\nabla u_h,\nabla v_h)=(f,v_h) \Forall v_h\in V_h.
\end{equation}
The n-c SPLS discretization \eqref{abstract:variational2_h} to be solved is: Find $(w_h,p_h=R_hA\nabla u_h)$ such that
\begin{equation}\label{discrete_problemSPLS}
\begin{array}{lclll}
(\nabla w_h,\nabla v_h) & + & (p_h,\nabla v_h) &= (f,v_h) &\ \Forall  v_h \in V_h,\\
R_hA\nabla w_h & & & =0.
\end{array}
\end{equation}
\subsubsection{Piecewise linear test space}\label{Vh=P1}

We make further assumptions to discuss stability for the family $\{(V_h,\M_h)\}$. We assume for simplicity $\Omega \subset \mathbb{R}^2$ is a polygonal domain separated into two subdomains by a smooth interface $\Gamma \subset \Omega$. The results can easily be extended to $N$ subdomains as well as polyhedral domains in $\mathbb{R}^3$. We also assume that the triangular mesh $\T_h$ is locally quasi-uniform. Let $\{z_{1,i},\dots,z_{N_i,i}\}$ be the set of all nodes of $\T_{h,i}$ and assume all triangles adjacent to $z_{j,i}$ are of regular shape and their area is of order $h_{j,i}^2$. In this notation, the mesh size of $\T_h=\T_{h,1}\cup \T_{h,2}$ is $h:=\max\{h_{1,1},h_{2,1},\dots, h_{N_1,1},h_{1,2},h_{2,2},\dots,h_{N_2,2}\}$. 

We take $V_h$ to be the space consisting of piecewise linear polynomials with respect to $\T_h$ vanishing on the boundary of $\Omega$. Also, we take $k=1$. Hence, each component of $M_{h,0}|_{\Omega_i}$ consists of continuous linear piecewise polynomials with respect to the mesh $\T_{h,i}$. Let $\{\Phi^i_1,...,\Phi^i_{2N_i}\}$ be a nodal basis for $M_{h,0}|_{\Omega_i}$ and assume that $\Phi^i_j=(\phi^i_j,0)^T$ and $\Phi^i_{N_i+j}=(0,\phi^i_j)^T$ for $j=1,\dots,N_i$. Here, $\{\phi^i_1,\dots,\phi^i_{N_i}\}$ is a nodal basis for the space of continuous piecewise linear polynomials with respect to $\T_{h,i}$. With this notation, we note that $\{A\Phi^1_j\}_{j=1}^{N_1} \cup \{A\Phi^2_j\}_{j=1}^{N_2}$ is a basis for $\tilde{\M}_h$. Lastly, we define $M_{A_i}$ to be the Gram matrix of the set $\{A\Phi^i_j\}_{j=1}^{N_i}$ with respect to the $(\cdot, \cdot)_{\tilde{Q}}$ inner product and $D_i:=\text{diag}\left(h_{1,i}^2,h_{2,i}^2,\dots, h_{N_i,i}^2,h_{1,i}^2,h_{2,i}^2,\dots, h_{N_i,i}^2\right)$. To prove stability for the family $\{(V_h,\M_h)\}$, we need the following two lemmata. The first lemma follows from a similar result (for no interfaces)  proved in \cite{BQ17} and, for completeness, is restated using the notation and assumptions from this section.
\begin{lemma}
Under the assumptions of Section \ref{Vh=P1}, we have that for $i=1,2$
\begin{equation}\label{MA estimate}
\langle M_{A_i}\gamma,\gamma\rangle_e \leq c\,a_{max}\langle D_i\gamma,\gamma\rangle_e \Forall \gamma \in \mathbb{R}^{2N_i}.
\end{equation}
Consequently,
\begin{equation}\label{MAinverse estimate}
\langle M^{-1}_{A_i}\gamma,\gamma\rangle_e \geq \frac{c}{a_{max}}\langle D^{-1}_i\gamma,\gamma\rangle_e \Forall \gamma \in \mathbb{R}^{2N_i}.
\end{equation}
\end{lemma}
We note that the constant $c$ in the above lemma is generic and does not depend on $h$. The next result shows that \eqref{eq:ProjCoerces} is satisfied for the representation operator $R_h$ defined in this section.
\begin{lemma}\label{coer_interface}
Under the assumptions of Section \ref{Vh=P1}, there exists a constant $c$, independent of $h$, such that
\begin{equation}\label{interfaceCoercivity}
\|R_h A \nabla v_h\|_h  \geq c \,  \frac{a_{min}}{a_{max}}  \ \|A \nabla v_h\|_{\tilde{Q}} \Forall  v_h \in V_h.
\end{equation}
\begin{proof}
First, note that $\{A\Phi^1_1,\dots,A\Phi^1_{2N_1}\}$ and $\{A\Phi^2_1,\dots,A\Phi^2_{2N_2}\}$ are nodal bases for $\tilde{\M}_h|_{\Omega_1}$ and $\tilde{\M}_h|_{\Omega_2}$, respectively. Define $v_h^i:=v_h|_{\Omega_i}$ for $v_h\in V_h$. For a fixed $A \nabla v_h$ with $v_h \in V_h$ we define the dual vectors $G^1_h \in \R^{2N_1},G^2_h \in \mathbb{\R}^{2N_2}$ by
\begin{align*}
&(G^1_h)_i:= (A\nabla v^1_h, A \Phi^1_i)_{\tilde{Q}} = (A\nabla v^1_h, \Phi^1_i) \ i=1,...,2N_1,\\
&(G^2_h)_i:= (A\nabla v^2_h, A \Phi^2_i)_{\tilde{Q}} = (A\nabla v^2_h, \Phi^2_i)\ i=1,...,2N_2,
\end{align*}
and let 
\[
R_h A \nabla v_h =
\begin{cases}
      \sum\limits_{i=1}^{2N_1} \alpha_i A \Phi^1_i& \textrm{in} \ \ \Omega_1,\\\\
      \sum\limits_{i=1}^{2N_2} \beta_i A \Phi^2_i& \textrm{in} \ \ \Omega_2.\\
\end{cases} 
\]
Thus,  $\alpha =(\alpha_1, \alpha_2, \dots,\alpha_{2N_1})^T$ and $\beta=(\beta_1,\beta_2, \dots, \beta_{2N_2})^T$ are solutions to  
\[
M_{A_1} \: \alpha = G^1_h, \qquad \text{and} \qquad M_{A_2} \: \beta = G^2_h,  
\]
respectively. Using \eqref{MAinverse estimate}, we obtain
\begin{align*}
\|R_h A \nabla v_h\|^2_h  &= \displaystyle \sum_{i, j =1}^{2N_1} \alpha_i\: \alpha_j \left( A \Phi^1_i, \Phi^1_j\right)+ \displaystyle \sum_{i, j =1}^{2N_2} \beta_i\: \beta_j \left( A \Phi^2_i, \Phi^2_j\right)\\
&= \left< M_{A_1}^{-1} G^1_h, G^1_h \right>_e+ \left< M_{A_2}^{-1} G^2_h, G^2_h \right>_e\\
&\ge \frac{1}{a_{max}} \left< D_1^{-1} G^1_h, G^1_h \right>+\frac{1}{a_{max}} \left< D_2^{-1} G^2_h, G^2_h \right>\\
&= \frac{1}{a_{max}} \sum_{i = 1}^{2N_1} h_{i,1}^{-2} \left( A\nabla v^1_h,\Phi^1_i \right)^2+ \frac{1}{a_{max}} \sum_{i = 1}^{2N_2} h_{i,2}^{-2} \left( A\nabla v^2_h,\Phi^2_i \right)^2. \\
\end{align*}
We recall by definition of $D_1,D_2$ that we have $h_{i,1}=h_{i+N_1,1}$ for $i=1,\dots, N_1$ and $h_{i,2}=h_{i+N_2,2}$ for $i=1,\dots,N_2$ in the above. Note that  
\begin{align*}
\frac{1}{a_{max}} \sum_{i = 1}^{2N_1} h_{i,1}^{-2} \left( A\nabla v^1_h,\Phi^1_i \right)^2
&=\frac{1}{a_{max}}
		\sum_{i = 1}^{N_1} h_{i,1}^{-2} \left[ \left(
		 a_{11} \frac{\partial v^1_h}{\partial x} + a_{12} \frac{\partial v^1_h}{\partial y}, \phi^1_i \right)^2  + \left(
 		 a_{21} \frac{\partial v^1_h}{\partial x} + a_{22} \frac{\partial v^1_h}{\partial y}, \phi^1_i \right)^2 \right]\\
&=\frac{1}{a_{max}}  \sum_{i = 1}^{N_1} h_{i,1}^{-2}
    		\sum_{\tau^1\subset supp(\phi_i)} \Bigg| 
    		\left( \begin{array}{cc}
		      (a_{11}, \phi^1_i)_{\tau^1} & (a_{12}, \phi^1_i)_{\tau^1} \\
		      (a_{21}, \phi^1_i)_{\tau^1} & (a_{22}, \phi^1_i)_{\tau^1} \end{array} \right)
		\begin{pmatrix}
		      \frac{\partial v^1_h |_{\tau^1}}{\partial x}\\
	    	      \frac{\partial v^1_h|_{\tau^1}}{\partial y}
	          \end{pmatrix} \Bigg|^2 \\
&\geq c_1\frac{a^2_{min}}{a_{max}} \sum_{i = 1}^{N_1} \sum_{\tau^1\subset supp(\phi_i)} h_{i,1}^{2} |\nabla v^1_h|^2_{\tau^1}\\
&= c_1\frac{a^2_{min}}{a_{max}}||\nabla v^1_h||^2_{\Omega_1},
\end{align*}
where the inequality above follows as the lowest eigenvalue of the matrix
\[
\displaystyle \begin{pmatrix}
		      (a_{11}, \phi^1_i)_{\tau^1} & (a_{12}, \phi^1_i)_{\tau^1} \\
		      (a_{21}, \phi^1_i)_{\tau^1} & (a_{22}, \phi^1_i)_{\tau^1} \end{pmatrix}
\]  
is bounded below by $c_1\, h_{i,1}^2 a_{min}$ with a constant $c_1$ independent of $\tau^1$ and $h$. Similarly, we can show
\[
\frac{1}{a_{max}} \sum_{i = 1}^{2N_2} h_{i,2}^{-2} \left( A\nabla v^2_h,\Phi^2_i \right)^2 \geq c_2\frac{a^2_{min}}{a_{max}}||\nabla v^2_h||^2_{\Omega_2}.
\]
Thus, 
\begin{align*}
\|R_h A \nabla v_h\|^2_h&\geq \min(c_1,c_2) \frac{a^2_{min}}{a_{max}}(||\nabla v^1_h||^2_{\Omega_1}+||\nabla v^2_h||^2_{\Omega_2})\\
&\geq \min(c_1,c_2) \frac{a^2_{min}}{a_{max}}||\nabla v_h||^2\\
&\geq \min(c_1,c_2) \frac{a^2_{min}}{a^2_{max}}||A\nabla v_h||^2_{\tilde{Q}}.
\end{align*}
For the last inequality, we use that
\[
\Vert A \nabla v_h \Vert_{\tilde{Q}}^2 = \left( A \nabla v_h, \nabla v_h \right) \le  a_{max} \Vert \nabla v_h \Vert ^2.
\]
\end{proof}
\end{lemma}
As a consequence of Lemma \ref{coer_interface}, equation \eqref{infsup_discrete}, and Proposition \ref{discrete-Stability}, we have the following result.
\begin{theorem}
Let $\Omega \subset \mathbb{R}^2$ be a polygonal domain and $\{T_h\}$ be a family of locally quasi-uniform meshes for $\Omega$. For each $h$, let $V_h$ be the space of continuous linear functions with respect to the mesh $\{\T_h\}$ that vanish on $\partial \Omega$ and $\M_h$ be the corresponding projection type trial space defined in Section \ref{projForInt}. Then the family of spaces $\{(V_h,\M_h)\}$ is stable.
\end{theorem}

We note that in the case when $A=I$,  we have that {\bf Step 1} of our Uzawa type iterative process  (of Section \ref{sec:solvers}) coincides with a standard gradient recovery technique with projection operator $R_h$  for solving the Laplace equation. The n-c SPLS iterative process goes beyond the projection of {\bf Step 1}. By computing further  $p_j$ iterations in {\bf Step 2}, we approach $R_h \nabla u_h$, which according to Proposition \ref{approxProj}, is a quasi-optimal approximation of $\nabla u$ with functions in $\M_h:=R_h \nabla V_h$.   

\section{Numerical Results}\label{Sec:NR}
We implemented the n-c SPLS discretization on second order elliptic PDE of the form \eqref{interface_prob}. For all of the examples presented, we took $\Omega$ to be a bounded polygonal or polyhedral domain and chose the test space $V_h\subset H_0^1(\Omega)$ to be the space of continuous piecewise linear polynomials with respect to the quasi-uniform, or locally quasi-uniform, meshes $\T_h$. The trial spaces are of the {\it projection type} as presented in Section \ref{projForInt}, and Algorithm \ref{alg:UCG} was used for all examples.

Based on the first inequality of \eqref{alg_estimate}, we used a stopping criterion of 
\[
\|q_j\| \leq c_0 h^2,
\] 
on each level for the case of convex domains and uniform refinement. This is because the maximum possible order for the discretization error $\| A\nabla u -  R_h A \nabla u_h \|$ would be order two. In the case of non-uniform refinement, we use a stopping criterion of 
\[
\|q_j\| \leq c_0 N_{dof}^{-2},
\] 
on each level, where $N_{dof}$ is the number of degrees of freedom. 

In practice, we notice that we cannot achieve order two. This could be because on each subdomain we approximate, in a weighted $L^2$ norm, a possibly smooth component of the flux, but use {\it subspaces of $C_0-P1$ functions} as approximation spaces component wise.

\begin{remark}
We note that for the SPLS discretization of the interface problem, the  primal variable $u$ can be approximated along the process simultaneously by separately storing the $u_j$ part of the iterates $p_j=R_h (A \nabla u_j)$, which can serve as a proxy $p_j$, and follow the updates for $p_j$ as in the algorithm.  However, for the piecewise linear approximation we consider here, we do not observe a higher order of approximation for the primal variable. We obtain a convergence rate of order one in the energy norm.
\end{remark}

\subsection{Interface problems}
In all examples presented, the constant $c$ will denote the size of the jump in the coefficients of the matrix $A$. The level of mesh refinement will be denoted by $k$. 

\subsubsection{Intersecting interface example}
For the first example, we took $\Omega = (0, 1) \times (0, 1)$ with the interface $\Gamma := \Omega \cap \{(x, y)\,|\, \ x=1/2 \ \text{or} \ y=1/2\}$ as considered in \cite{Bank-Xu03}. The family of interface-fitted, locally quasi-uniform meshes $\{\T_h\}$ was obtained by a standard uniform refinement strategy starting with a uniform coarse mesh. We computed $f$ such that for
\[
 A(x,y) =a(x,y) I_{2}, \ \text{where}\ 
\begin{aligned}
 a(x,y) &= 
	\begin{cases} 
	 1 \  & \mbox{if} \; (x,y) \in [0,1/2]^2 \cup [1/2,1]^2,\\ 
	 c & \mbox{if} \; (x,y) \in \Omega \setminus ( [0,1/2]^2 \cup [1/2,1]^2),\\ 
	\end{cases} \\
\end{aligned}
\] 
the exact solution is $u(x,y)=a(x,y)^{-1}\sin(2\pi x)\sin(2\pi y)$. Table \ref{table_interface1} shows the results for $c=1/10,1/100,$ and $1/1000$.
\begin{table}[!ht!]
\begin{center}
$\mathrm{error}=\|A\nabla u-R_hA\nabla u_h\|_{\tilde{Q}}$
\resizebox{12.5cm}{1.7cm}{
\begin{tabular}{|*{10}{c|}}
\hline
\multicolumn{1}{|c|}{ \multirow{2}{*}{\parbox{1.4cm}{\centering $h = 2^{-k}$ $k$}} } & \multicolumn{3}{c|}{$c=1/10$} & \multicolumn{3}{c|}{$c=1/100$}
& \multicolumn{3}{c|}{$c=1/1000$} \\  
\cline{2-10}
 & error &  rate  & it & error & rate & it & error & rate  & it \\
\hline
1 & 5.177 &           &4 & 15.686 &             & 4 & 49.383 &            & 4 \\
\hline
2 & 1.262 & 2.037 & 10 & 3.947 & 1.990 & 12 & 15.827 & 1.642 & 11 \\
\hline
3 & 0.339 & 1.895 & 16 & 1.070 & 1.882 & 27 & 3.607 &  2.134 & 29 \\
\hline
4 & 0.097 & 1.802 & 17 & 0.307 & 1.803 & 33 & 0.985 & 1.873 & 63 \\
\hline
5 & 0.027 & 1.849 & 22 & 0.086 & 1.832 & 44 & 0.295 & 1.738 & 76 \\
\hline
\end{tabular}
\caption{Interface problem with intersecting interfaces.}\label{table_interface1}
}
\end{center}
\end{table}

\subsubsection{Gradient singularity at the origin}
For the second example, we solved \eqref{interface_prob} where the gradient of the solution is singular at the origin, see \cite{petzoldt2001regularity}. The domain $\Omega = (-1, 1)^2$ is decomposed as $\Omega_2:=\{(x,y)\in \Omega \,|\, 0 <\theta(x,y) < \pi/2\}$ and $\Omega_1:=\Omega \setminus \Omega_2$, where $\theta(x,y)$ is the angle in polar coordinates of the point $(x,y)$. We computed $f$ such that for
\[
 A(x,y) =a(x,y) I_{2}, \ \text{where}\ 
\begin{aligned}
 a(x,y) &= 
	\begin{cases} 
	 1 \  & \mbox{if} \; (x,y) \in \Omega_1,\\ 
	 c & \mbox{if} \; (x,y) \in \Omega_2,\\ 
	\end{cases} \\
\end{aligned}
\] 
the exact solution, given in polar coordinates, is $u(r,\theta)= r^{\lambda}(1-r)^2\mu(\theta)$ where 
\[
\begin{aligned}
 \mu(\theta) &= 
	\begin{cases} 
	 \cos(\lambda(\theta-\pi/4)) \  & \mbox{if} \; (x,y) \in \Omega_2, \\ 
	 b\cos(\lambda(\pi-|\theta-\pi/4|)) \ & \mbox{otherwise},\\ 
	\end{cases} \\
\end{aligned}
\] 
and 
\[
\lambda = \frac{4}{\pi}\arctan\left(\sqrt{\frac{3+c}{1+3c}}\right), \ \ \ \ b=-c\frac{\sin\left(\lambda \frac{\pi}{4}\right)}{\sin\left(\lambda \frac{3\pi}{4}\right)}.
\]

Using a similar standard uniform refinement strategy as in the previous problem, Table \ref{table_interface2u} summarizes results for $c=5$ and $c=15$.

\begin{table}[!ht!]
\begin{center}
$\mathrm{error}=\|A\nabla u-R_hA\nabla u_h\|_{\tilde{Q}}$
\resizebox{9.5cm}{1.7cm}{
\begin{tabular}{|*{10}{c|}}
\hline
\multicolumn{1}{|c|}{ \multirow{2}{*}{\parbox{1.4cm}{\centering $\mathrm{level} \ k$}}} & \multicolumn{3}{c|}{$c=5$} & \multicolumn{3}{c|}{$c=15$}\\  
\cline{2-7}
 & error &  rate  & it & error & rate & it \\
\hline
1 & 1.230 &           &4     & 3.397 &             & 5 \\
\hline
2 & 0.372 & 1.726 & 10 & 1.199 & 1.503 & 16 \\
\hline
3 & 0.180 & 1.049 & 17 & 0.662 & 0.856 & 44 \\
\hline
4 & 0.101 & 0.837 & 36 & 0.402 & 0.721 & 104 \\
\hline
5 & 0.058 & 0.795 & 57 & 0.246 & 0.706 & 166 \\
\hline
\end{tabular}
\caption{Interface problem with gradient singularity at $(0,0)$ on uniform mesh.}\label{table_interface2u}
}
\end{center}
\end{table}

Using uniform meshes, we observe a convergence rate less than one. To better capture the singularity of the gradient, a family of interface-fitted, locally quasi-uniform meshes $\{\T_h\}$ was obtained by a graded refinement strategy depending on a refinement parameter $\kappa$ \cite{BNZ1, BNZ2}. The refinement is done by splitting each triangle in four smaller triangles. In particular, we divide every edge that contains the singular point (the origin in this case) under a fixed ratio $\kappa$ such that the edge containing the singular point is $\kappa$ times the other segment. In the case $\kappa=1$, we recover the uniform refinement. Numerical results using graded meshes with $\kappa = 0.22$ are summarized in Table \ref{table_interface2} for $c=5$ and $c=15$.

\begin{figure}
\includegraphics[scale=.3]{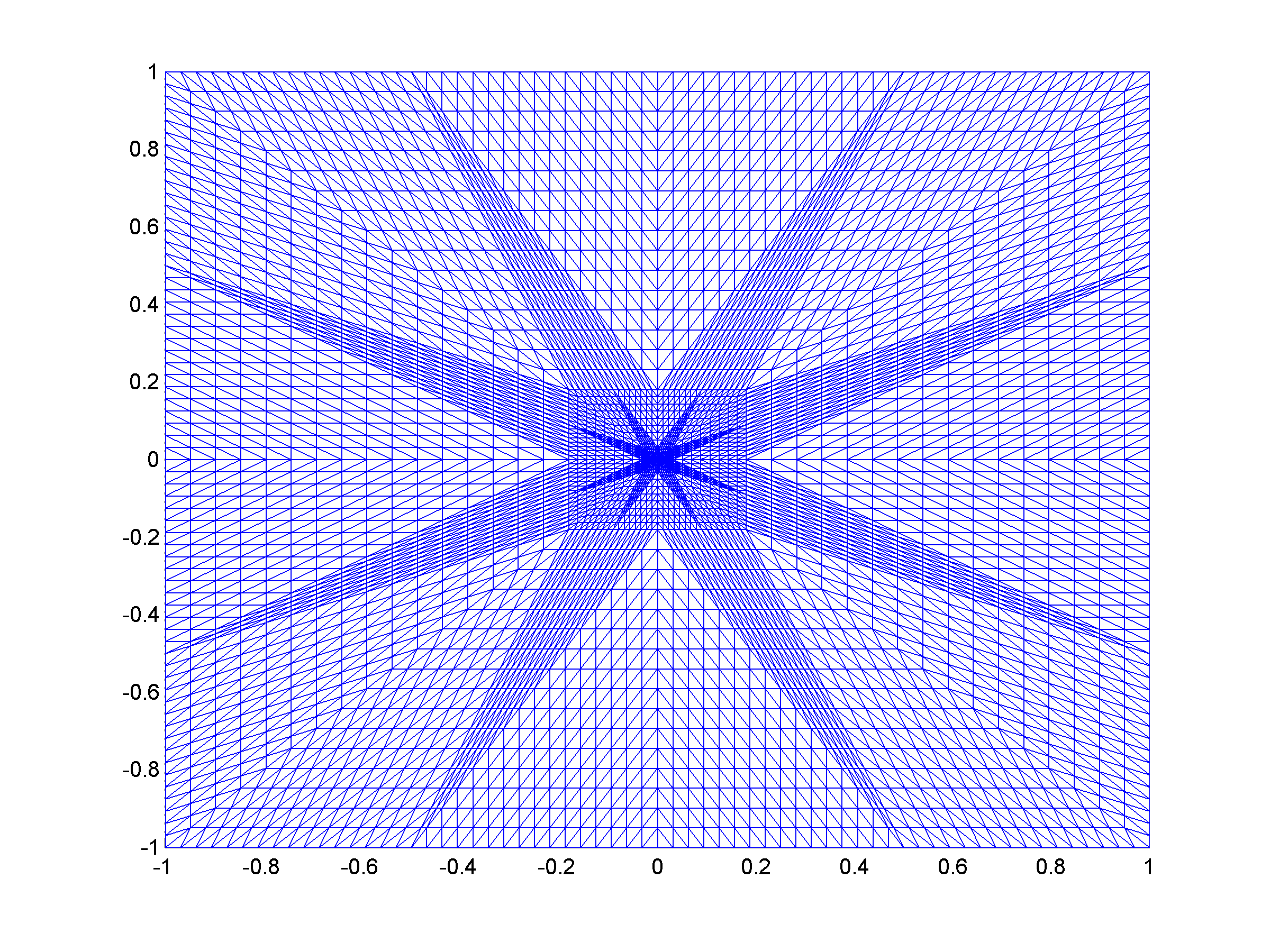} 
\includegraphics[scale=.3]{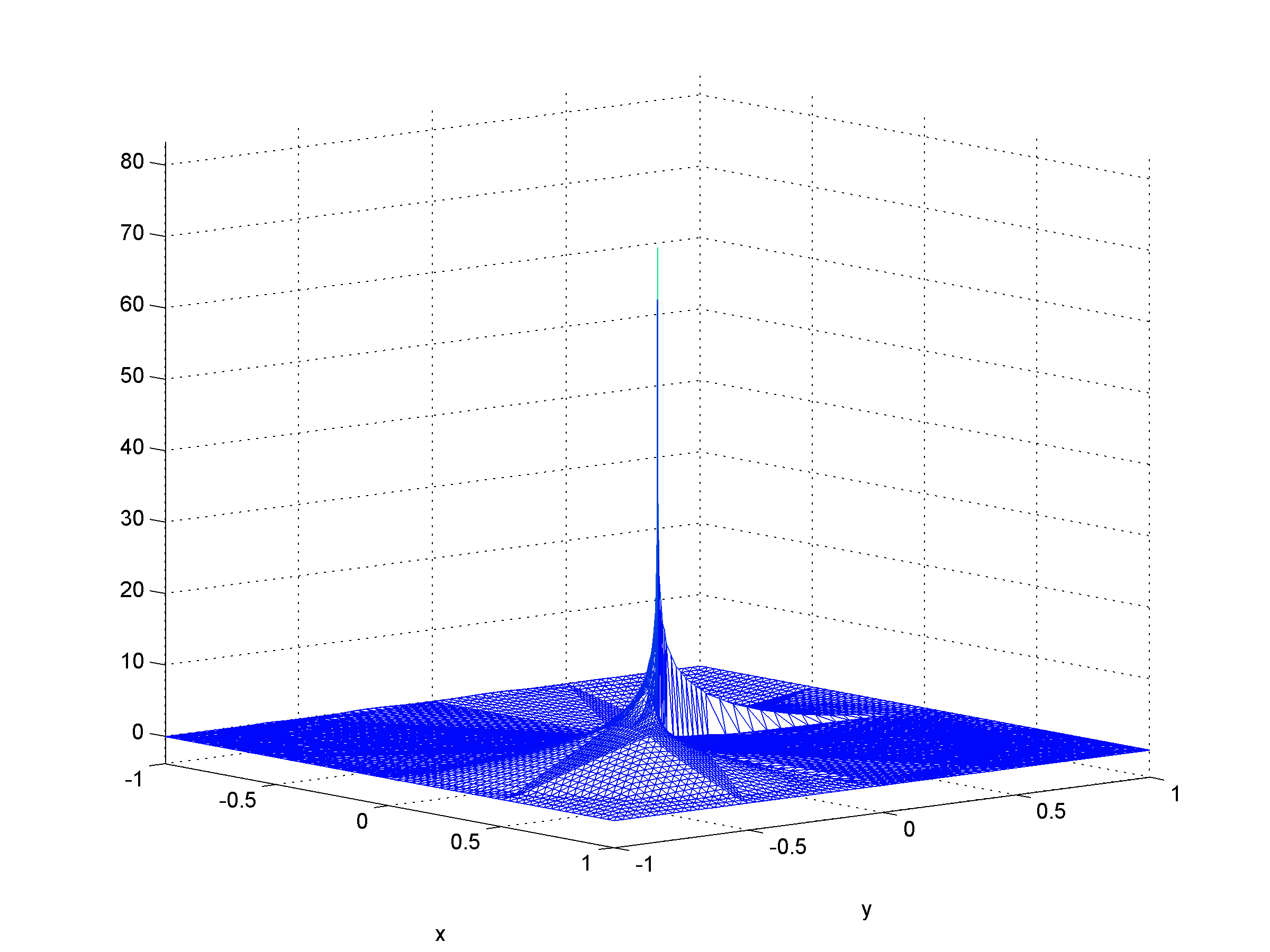} 
\caption{Mesh and $x$ component of the computed flux for gradient singularity example.}
\label{fig_interface}
\end{figure}

\begin{table}[!ht!]
\begin{center}
$\mathrm{error}=\|A\nabla u-R_hA\nabla u_h\|_{\tilde{Q}}$
\resizebox{9.5cm}{1.7cm}{
\begin{tabular}{|*{10}{c|}}
\hline
\multicolumn{1}{|c|}{ \multirow{2}{*}{\parbox{1.4cm}{\centering $\mathrm{level} \ k$}}} & \multicolumn{3}{c|}{$c=5$} & \multicolumn{3}{c|}{$c=15$}\\  
\cline{2-7}
 & error &  rate  & it & error & rate & it \\
\hline
1 & 0.949 &            &4 & 2.605 &           & 5 \\
\hline
2 & 0.585 & 0.699 & 9 & 1.504 & 0.792 & 15 \\
\hline
3 & 0.151 & 1.945 & 16 & 0.412 & 1.868 & 46 \\
\hline
4 & 0.052 & 1.545 & 23 & 0.143 & 1.529 & 72 \\
\hline
5 & 0.017 & 1.602 & 31 & 0.047 & 1.600 & 94 \\
\hline
\end{tabular}
\caption{Interface problem with gradient singularity at $(0,0)$ on non-uniform mesh.}\label{table_interface2}
}
\end{center}
\end{table}
Figure \ref{fig_interface} depicts the mesh generated (at the final level of refinement) using the graded refinement strategy for $\kappa=0.22$ as well as the $x$ component of the computed gradient for the case of $c=15$.

\subsubsection{3-$D$ example}
For the third example, we took $\Omega \subset \mathbb{R}^3$ to be the unit cube with interface $\Gamma:=\Omega \cap \{(x,y,z)\,|\,x=1/2\}$. We computed $f$ such that for 
\[
 A(x,y,z) =a(x,y,z) I_{3}, \ \text{where}\ 
\begin{aligned}
 a(x,y,z) &= 
	\begin{cases} 
	 1 \  & \mbox{if} \; x < \frac{1}{2},\\ 
	 c & \mbox{if} \; x \geq \frac{1}{2},\\ 
	\end{cases} \\
\end{aligned}
\] 
the exact solution is 
\[
\begin{aligned}
 u(x,y,z) &= 
	\begin{cases} 
	 c \,  x(x-\frac{1}{2})y(y-1)z(z-1) & \mbox{if} \; x < \frac{1}{2}, \\ 
	(x-\frac{1}{2})(x-1)y(y-1)z(1-z) & \mbox{if} \; x \ge \frac{1}{2}.\\ 
	\end{cases} \\
\end{aligned}
\]
Table \ref{table_interface3} shows the results for $c=5,25,$ and $50$.
\begin{table}[!ht!]
\begin{center}
$\mathrm{error}=\|A\nabla u-R_hA\nabla u_h\|_{\tilde{Q}}$
\resizebox{12.5cm}{1.7cm}{
\begin{tabular}{|*{10}{c|}}
\hline
\multicolumn{1}{|c|}{ \multirow{2}{*}{\parbox{1.4cm}{\centering $h = 2^{-k}$ $k$}} } & \multicolumn{3}{c|}{$c=5$} & \multicolumn{3}{c|}{$c=25$}
& \multicolumn{3}{c|}{$c=50$} \\  
\cline{2-10}
 & error &  rate  & it & error & rate & it & error & rate  & it \\
\hline
1 & 0.0456 &           &1 & 0.2124 &           & 1 & 0.4208 &             & 1 \\
\hline
2 & 0.0159& 1.517 & 6 & 0.0744& 1.513 & 14 & 0.1475 & 1.512 & 18 \\
\hline
3 & 0.0042 & 1.925 & 9 & 0.0196 & 1.922 & 27 & 0.0389 &  1.921 & 44 \\
\hline
4 & 0.0011 & 1.879 & 12 & 0.0053 & 1.882 & 42 & 0.0106 & 1.881 & 67 \\
\hline
5 & 0.0003 & 1.863 & 15 & 0.0014 & 1.889 & 65 & 0.0028 & 1.889 & 110 \\
\hline
\end{tabular}
\caption{3-$D$ interface problem.}\label{table_interface3}
}
\end{center}
\end{table}

We observe for  both convex and non-convex  domains that the approximation of the flux is {\it super-linear}, and the method works well no matter the size of the jump discontinuity. Also, we notice that the number of iterations depends on the size of the jump as well as $h$. This is in accordance with Remark \ref{U-remark} in that the number of iterations on each level will be proportional with $\frac{M}{m_h}$, which depends on the jump $\frac{a_{max}}{a_{min}}$  and $-log(h^2)$. 

\subsection{Flux recovery for highly oscillatory coefficients} 
We note that the stability and approximation results of Section \ref{sec:Mh} can be applied to the case when the coefficients 
of the PDE (the entries of $A$) are smooth functions. We would like to illustrate the advantages of our n-c SPLS discretization with projection on an example where the matrix $A$ has highly oscillatory coefficients.  We solved \eqref{interface_prob} on $\Omega = (0, 1) \times (0, 1)$ with $A=a(x,y)I_2$, where 
\[
a(x,y)=\frac{1}{4+P(\sin(2\pi x/\varepsilon)+\sin(2\pi y/\varepsilon))}.
\]
We computed $f$ such that the exact solution is given by
\[
u(x,y)=\frac{\sqrt{4-P^2}}{2}(x^2+y^2)\exp\left(\frac{1}{x^3-x}+\frac{1}{y^3-y}\right).
\]
This is a small modification of a similar  example presented in \cite{Mu-Wang-Li16}. Table \ref{hoc_example} shows the results for various values of $\varepsilon$. In all computations, we chose $P=1.8$.
\begin{table}[!ht!]
\begin{center}
$\mathrm{error}=\|A\nabla u-R_hA\nabla u_h\|_{\tilde{Q}}$
\resizebox{12.5cm}{1.7cm}{
\begin{tabular}{|*{10}{c|}}
\hline
\multicolumn{1}{|c|}{ \multirow{2}{*}{\parbox{1.4cm}{\centering $h = 2^{-k}$ $k$}} } & \multicolumn{3}{c|}{$\varepsilon=0.2$} & \multicolumn{3}{c|}{$\varepsilon=0.1$} & \multicolumn{3}{c|}{$\varepsilon=0.05$}\\  
\cline{2-10}
 & error & rate  & it & error & rate  & it & error & rate  & it\\
\hline
5 & 1.34e-04 & 2.96 & 4 &1.71e-04 & 1.88 & 3 & 3.23e-04 & 0.27 & 2\\
\hline
6 & 5.65e-05 & 1.24 & 6 & 5.46e-05 & 1.65 & 5 & 6.44e-05 & 2.33 & 4\\
\hline
7 & 1.42e-05 &  1.99 & 9 & 1.34e-05 & 2.03 & 8& 1.22e-05 & 2.40 & 7\\
\hline
8 & 4.07e-06 & 1.80 & 12 & 2.57e-06 & 2.39 & 12 & 2.38e-06 & 2.36 & 11\\
\hline
\end{tabular}
\caption{Highly oscillatory coefficients example.}\label{hoc_example}
}
\end{center}
\end{table}

The numerical results show almost $O(h^2)$ order of approximation for the flux for meshes that are small enough to capture the high frequency of the coefficients due to the size of  $\varepsilon$. 

\section{Conclusion}
We presented a  saddle point least squares method with non-conforming trial spaces for discretization of second order PDEs with discontinuous coefficients. The proposed method is easy to implement using Uzawa type algorithm and leads to higher order approximation of the flux  if compared with standard finite element (non-mixed) techniques  based on  linear element approximation. In addition, the method works well when solving second order problems with variable coefficients, including highly oscillatory coefficients, and can be combined with known gradient recovery techniques and graded meshes techniques in order to construct optimal or quasi-optimal  discrete approximation spaces for the flux. 

We plan to further combine the n-c SPLS method with known multilevel and adaptive techniques \cite{ainsworth-oden, carstensen-hoppe-lobhardUnifAdapt12} for  designing  robust iterative solvers  for more general second order elliptic PDE that are parameter dependent and exhibit singular solutions due to non-convex domains or discontinuous coefficients.


\bibliography{../../bacutaBib}
\bibliographystyle{plain} 

\end{document}